\documentclass[10pt]{amsart}
\usepackage{palatino, amssymb, amsfonts, latexsym}
\usepackage{amssymb, amsmath, amscd}
\input xy
\xyoption{all}

\newtheorem{theorem}{Theorem}[section]
\newtheorem{lemma}[theorem]{Lemma}
\newtheorem{prop}[theorem]{Proposition}
\newtheorem{cor}[theorem]{Corollary}
\theoremstyle{definition}
\newtheorem{definition}[theorem]{Definition}

\theoremstyle{remark}
\newtheorem{remark}[theorem]{Remark}

\newcommand{\U}{\mathbf{U}}

\newcommand{\A}{\mathcal A}

\newcommand{\Ud}{\dot{\mathbf U}}
\newcommand{\Uda}{\dot{\mathbb U}}
\newcommand{\B}{\mathbf B}
\newcommand{\Bd}{\dot{\mathbf B}}

\numberwithin{equation}{section}

\begin{document}

\title{Generalized $q$-Schur algebras and quantum Frobenius.}

\author{Kevin McGerty}
\address{Department of Mathematics, University of Chicago. }

\date{August, 2006}

\begin{abstract}
The quantum Frobenius map and it splitting are shown to descend to maps between generalized $q$-Schur algebras at a root of unity. We also define analogs of $q$-Schur algebras for any affine algebra, and prove the corresponding results for them. 
\end{abstract}

\maketitle

\section{Introduction}
In their study of the quantum group $\mathbf U(\mathfrak{sl}_n)$, Beilinson, Lusztig and MacPherson defined a family of finite-dimesional quotients known as $q$-Schur algebras (see also \cite{DJ}). Subsequently Lusztig defined similar families of finite dimensional quotients for any quantum group of finite type, which here, following \cite{Doty}, we call ``generalized $q$-Schur algebras''. We show that, at a root of unity, the quantum Frobenius map and a splitting of it (which extends the splitting map defined by Lusztig) descend to a family of maps between generalized $q$-Schur algebras (Theorem $4.5$ below). Note that in the general situation, the quantum Frobenius map is a homomorphism from the quantum group at a root of unity to the integral enveloping algebra of a different quantum group (the Langlands' dual in the case where the order of the root of unity is divisible by $6$ for example). For the original $q$-Schur algebras we also give a simple characterization of the embedding in the geometric context of \cite{BLM90} (Lemma 5.2). Finally we define a family of algebras which are affine analogs of generalized $q$-Schur algebras (essentially an algebraic version of the algebras considered in the work of Nakajima), and establish the analogous results for these algebras (Theorem 6.2).

One impetus for this work came from discovery by Fayers and Martin \cite{FM} of an explicit embedding of the Schur algebra $S_{\mathbf k}(2,r)$ into the Schur algebra $S_{\mathbf k}(2,pr)$ when $\mathbf k$ is a field of characteristic $p$. This paper grew out of a desire to find a quantum version of their embedding, which turns out to be provided by the splitting map for the quantum Frobenius (see Section \ref{n=2case}). Our results show that the quantum Frobenius is compatible with certain families of based ideals in the modified quantum group, suggesting that there are compatibilities between this map and the canonical basis that have yet to be understood. We hope to return to this issue in future work (also related to \cite{McG06}).

\section{Quantum Groups and Modified Quantum Groups}
\label{modified}
\subsection{Cartan and Root Data}
In this section we recall the definition of a quantum group and its associated modified form, using Lusztig's notion of a Cartan datum.
 
\begin{definition}
A \textit{Cartan datum} \cite[1.1.1]{L93} is a pair $(I,\cdot)$ consisting of a finite set
$I$ and a $\mathbb Z$-valued symmetric bilinear pairing on the
free Abelian group $\mathbb Z[I]$, such that
\begin{itemize}
\item $i\cdot i \in \{2,4,6,\ldots\}$
\item $2\frac{i\cdot j}{i\cdot i} \in \{0,-1,-2,\ldots\}$, for $i\neq j$.
\end{itemize}
We will write $a_{ij} = 2\frac{i\cdot j}{i \cdot i}$. Note that the matrix $A= (a_{ij})$ is a symmetrizable generalized Cartan matrix.
A \textit{root datum} \cite[2.2.1]{L93} of type $(I,\cdot)$ is a pair $Y,X$ of finitely-generated
 free Abelian groups and a perfect pairing
$\langle,\rangle \colon Y \times X \to \mathbb Z$, together with
imbeddings $I\subset X$, ($i\mapsto \alpha_i$) and $I\subset Y$, ($i
\mapsto \check{\alpha}_i$) such that $\langle \check{\alpha}_i,\alpha_j \rangle = 2\frac{i\cdot j}{i
\cdot i}$. Given a Cartan datum, one attach to it a Coxeter system $(W,S)$ where $S = \{s_i: i \in I\}$ is the set of simple reflections satisfying braid relations of length $h(i,j)$ where $\cos^2\frac{\pi}{h(i,j)} = \frac{i\cdot j}{i \cdot i}\frac{j\cdot i}{j \cdot j}$. The group $W$ (the Weyl group) acts on $X$ via $s_i\colon X \to X$ where $s_i(\lambda) = \lambda - \langle \check{\alpha}_i, \lambda \rangle \alpha_i$.
\end{definition}

When there is no danger of ambiguity we shall, by abuse of notation, write $(X,Y)$ for a root datum. 

\subsection{Quantum Groups}
Given a root datum, we define an associated quantum group $\mathbf U$. Let $v$ be an indeterminate, and for each $i \in I$ set $v_i = v^{(i\cdot i)/2}$. We set 
\[
[n]_i = (v_i^n - v_i^{-n})/(v_i-v_i^{-1}) \in \mathbb Z[v,v^{-1}],
\] 
and then define 
\[
[n]_i! = [n]_i[n-1]_i.\ldots[1]_i; \qquad {n \brack k}_i = \frac{[n]_i!}{[k]_i! [n-k]_i!},
\]
(note that in fact ${n \brack k}_i \in \mathbb Z[v,v^{-1}]$). Let $\mathbf U$ be the $\mathbb Q(v)$-algebra generated by symbols $E_i, F_i, K_\mu$, $i \in I$, $\mu \in Y$, subject to the following relations.
\begin{enumerate}
\item $K_0=1$, $K_{\mu_1}K_{\mu_2} = K_{\mu_1+\mu_2}$ for $\mu_1,\mu_2 \in Y$;
\item $K_{\mu} E_i K_{\mu}^{-1} = v^{\langle\mu,\alpha_i\rangle}E_i, \quad K_{\mu} F_i K_{\mu}^{-1} =
 v^{-\langle\mu,\alpha_i \rangle}F_i$ for all $i \in I$, $\mu \in Y$;
\item $E_iF_j - F_jE_i = \delta_{i,j}\frac{K_i-K_i^{-1}}{v_i-v_i^{-1}}$;
\item $\sum_{r+s=1-a_{ij}}(-1)^{r} E_i^{(r)}E_jE_i^{(s)} = 0 \text{ for } i \neq j$;
\item $\sum_{r+s=1-a_{ij}} (-1)^rF_i^{(r)}F_jF_i^{(s)} = 0 \text{ for } i \neq j.$
\end{enumerate}
where $E_i^{(r)}$ denotes the $q$-divided power $\frac{E_i^r}{[r]_i!}$, and $K_i = K_{\frac{1}{2}(i \cdot i) \check{\alpha}_i}$.

\subsection{Modified Quantum Groups}
For our construction it is better to work with a variant of $\U$ introduced by Lusztig, which we now describe (following \cite{K94}). Let $\mathbf{Mod}_X$ denote the category of left $\mathbf U$-modules endowed with a weight decomposition: the objects of $\mathbf{Mod}_X$ are $\U$-modules $V$ such that
$$ V = \bigoplus_{\lambda \in X} V_\lambda,$$
where
$$ V_\lambda = \{v \in V \colon K_\mu v = v^{\langle \mu,\lambda\rangle}v, \forall \mu \in Y\}.$$
Let $R$ be the endomorphism ring of the forgetful functor to the category of vector spaces. By definition, an element of $a$ of $R$ associates to each object $V$ of $\mathbf{Mod}_X)$ a linear map $a_V$, such that $a_W\circ f = f \circ a_V$ for any morphism $f\colon V \to W$.
Any element of $\mathbf U$ clearly determines an element of $R$. For each $\lambda \in X$,
let $1_\lambda \in R$ be the projection to the $\lambda$ weight space. 
\begin{definition}
Let $\Ud$ be the subring (in fact, clearly, $\mathbb Q(v)$-subalgebra)
$$\dot{\mathbf U} = \bigoplus_{\lambda \in X} \mathbf U 1_\lambda.$$
Note that $\Ud$ does not have a multiplicative identity, but instead a collection of orthogonal idempotents. It is clear that the category $\mathbf{Mod}_X$ is equivalent to a category of modules for $\Ud$, the category of \textit{unital} modules.  
\end{definition}

Let $\U^\pm$ be the (isomorphic) subalgebras of $\U$ generated by $\{E_i: i \in I\}$ and $\{F_i: i \in i\}$ respectively. It is known that if $\mathbf f$ is the associative $\mathbb Q(v)$-algebra generated by $\{\theta_i, i \in I\}$ subject only to the analog of relation $(4)$ for $\U$:
\[
\sum_{r+s=1-a_{ij}}(-1)^r \theta_i^{(r)}\theta_j\theta_i^{(s)} = 0 \text{ for } i \neq j,
\]
then the natural maps $\mathbf f \to \U^\pm$ are isomorphisms. 
For $x \in \mathbf f$, let $x^{\pm}$ denote the image of $x$ in $\mathbf U^{\pm}$. 

\subsection{Integral Forms}
Since we are interested in the algebra $\Ud$ at a root of unity, we need to specify an integral form of the algebra explicitly. There are two distinct integral forms for a quantum group -- with or without divided powers. In what follows we will only need the modified quantum group in its integral form with divided powers. Let $\A$ denote the ring $\mathbb Z[v,v^{-1}]$ of Laurent polynomials in $v$. By \cite[23.2]{L93} the $\A$-subalgebra $\Ud_\A$ of $\Ud$ generated by $E_i^{(n)}1_\lambda, F_i^{(n)}1_\lambda$ $(i \in I, n \geq 0, \lambda \in X$) is an integral form (\textit{i.e.} the canonical map $\mathbb Q(v) \otimes_{\A} \Ud_\A \to \Ud$ is an isomorphism).  Both the algebra $\Ud_\A$ and the corresponding integral form $\mathbf f_\mathcal A$ of $\mathbf f$ have canonical bases $\Bd$ and $\B$ (see \cite[25.2]{L93} and \cite[14.4]{L93} respectively).

Ley $\ell$ be any positive integer. Following \cite[35.1.3]{L93}, for $\ell$ even we set  $l=2\ell$, while if $\ell$ is odd we set $l= \ell$ or $2\ell$. We set $\A_\ell$ to be the quotient ring $\A/(\Phi_{l}(v))$ where $\Phi_l$ is the $l$-th cyclotomic polynomial, and then let $\Ud_{\ell}$ be the corresponding specialization $\A_\ell \otimes_{\A} \Ud_{\A}$ of the modified form $\Ud$. More generally for any $\A$-algebra $R$ we will write $_R\Ud$ for the corresponding specialization of $\Ud_\A$, and $_R\mathbf f$ for the specialization of $\mathbf f_\A$.

In order for us to define the embeddings we must also introduce the notion of an $\ell$-modified root datum \cite[2.2.4]{L93}:

\begin{definition} Given a Cartan datum $(I,\cdot)$ and an integer $\ell$ we define a new Cartan datum $(I,\circ)$ given by
\[
i \circ j = l_il_j(i\cdot j),
\]
where $l_i$ is the smallest positive integer such that $l_i(i\cdot i/2) \in \ell \mathbb Z$. It is easy to check that $(I, \circ)$ is indeed a Cartan datum. Given a root datum $(X,Y)$ we can also define a new root datum by setting $X^*= \{\lambda \in X\colon \langle \check{\alpha}_i, \lambda \rangle \in l_i \mathbb Z\}$ and $Y^* = \text{Hom}(X^*, \mathbb Z)$, with the obvious pairing between $X^*$ and $Y^*$. The simple roots of the $(X^*,Y^*)$ are $\alpha_i^* = l_i \alpha_i$ and the simple coroots are $\check{\alpha}_i^*$, where $\check{\alpha}_i^*(\lambda) = l_i^{-1}\langle \check{\alpha}_i, \lambda \rangle$. Note that there is thus a natural inclusion $X^* \to X$ and hence a restriction map $Y \to Y^*$.
\end{definition}

In the case where the root datum is of finite type (\textit{i.e.} the matrix $A$ is positive definite), the integers $i\cdot i/2$ lie in  $\{1,2,3\}$ and so for $\ell$ divisible by $6$, the root system $(Y^*, X^*)$ is essentially that attached to the dual group. If $\Ud$ is the modified quantum group attached to a root datum $(X,Y)$ we denote $v_i^* = v^{(i \circ i)/2} = v_i^{l_i^2}$, and the modified quantum group attached to the root datum $ (X^*,Y^*)$ by $\Ud^*$. For clarity we will also denote the generator of $\Ud^*$ by $e_i^{(n)}1_\lambda$, $f_i^{(n)}1_\lambda$. Since $(v_i^*)^2 = (v_i^{l_i^2})^2 = 1$ in $\A_\ell$ the algebra $\Ud^*_\ell$ is close to the classical enveloping algebra. In \cite[33.2]{L93} this is called the \textit{quasiclassical case}. It is shown there \cite[Proposition 33.2.3]{L93} that a quasiclassical specialization of a quantum group satisfying mild hypotheses (satisfied for example by any finite type quantum group) is in fact isomorphic to the classical form $\Ud_{1}$, \textit{i.e.} where $v \mapsto 1$, thus the quantum Frobenius can be regarded as a map from $\Ud_\ell$ to $\mathcal A_\ell \otimes_\mathbb Z \Ud^{*}_1$.

\section{The contracting homomorphism}
\label{contract}

\subsection{The Quantum Frobenius}
We wish to construct embeddings of generalized $q$-Schur algebras, which are quotients of the modified quantum group. In this section we construct, as a first step, an embedding of the algebra $\Ud^*_\ell$ into $\Ud_\ell$. For this we rely on the work of Lusztig on the quantum Frobenius homomorphism. Recall from \cite[chapter 35]{L93} that there are two $\A_\ell$-homomorphisms $Fr\colon \mathbf f_{\mathcal A_\ell} \to \mathbf f^*_{\mathcal A_\ell}$ and $Fr' \colon \mathbf f^*_{\mathcal A_\ell} \to \mathbf f_{\mathcal A_\ell}$ which are given on generators by $Fr'(\theta_i^{(n)}) = \theta_i^{(nl_i)}$, and $Fr(\theta_i^{(n)}) = \theta_i^{(n/l_i)}$ when $l_i$ divides $n$ and zero otherwise. Lusztig also shows that $Fr$ ``extends'' to a map $Fr\colon \Ud_\ell \to \Ud^*_\ell$ where $Fr(E_i^{(n)}1_\lambda) = e_i^{(n/l_i)}1_\lambda$, if $l_i$ divides $n$ and $\lambda \in X^*$, and to zero otherwise, and similarly $Fr(F_i^{(n)}1_\lambda) = f_i^{(n/l_i)}1_\lambda$, if $l_i$ divides $n$ and $\lambda \in X^*$, and to zero otherwise. Note that $Fr$ is obviously surjective.  The observation of this section (which appears to be new) is that the map $Fr'$ also has an ``extension'' to $\Ud^*_\ell$.

\begin{remark}
The map $Fr$ is a $q$-analogue of the Frobenius morphism on enveloping algebras in positive characteristic, while the map $Fr'$ should be thought of as a $q$-analogue of the canonical ``Frobenius splitting'' map for Schubert varieties \cite{KL}.
\end{remark}

The following two lemmas are discussed in \cite{L93} and similar calculations appear in \cite{Li}.

\begin{lemma}
\label{binom}
Let $\phi\colon \A \to R$ be an algebra over $\A$ such that $\phi(v^{2l}) =1$, but $\phi(v^{2t}) \neq 1$ for all $t < \ell$. Then in $R$ we have 
\begin{enumerate}
\item ${m \brack k}= 0$ when $m$ is a multiple of $\ell$ and $k$ is not divisible by $\ell$.
\item Suppose that for $m \geq k$ we have $m = m_1\ell + s$ and $k = k_1\ell + t$ where $0 \leq s,t < \ell $.  Then we have
\[
{m \brack k} = v^{\ell(k_1s - m_1 t)+(m_1+1)k_1\ell^2}{m_1 \choose k_1}{s \brack t}
\]

\end{enumerate}
\end{lemma}
\begin{proof}
See Lemma $34.1.2$ parts $(a)$ and $(c)$ in \cite{L93}.

\end{proof}

\begin{lemma}
\label{present}
For any $R$-algebra $\phi\colon R \to \A$, the algebra $_R \Ud$ has the following presentation  (where we implicitly use the isomorphisms $_R \mathbf f \cong$ $_R\mathbf U^{\pm}$). Let $\zeta \in X$, then $_R \Ud$ is generated by elements $u^+1_\zeta u^-$ and $u^- 1_\zeta u^+$ for $u^\pm \in$ $_R\mathbf U^{\pm}$ such that:
\begin{itemize}
\item $E_i^{(a)} 1_\zeta F_j^{(b)} = F_j^{(b)}1_{\zeta + a\alpha_i + b\alpha_j} E_i^{(a)}$ for $i \neq j$;

\item $E_i^{(a)} 1_{-\zeta} F_j^{(b)} = \sum_{t \geq 0} \phi({a+b-\langle \check{\alpha}_i, \zeta\rangle \brack t}_i) F_i^{(b-t)} 1_{-\zeta +(a+b-t)\alpha_i}E_i^{(a-t)}$;

\item $F_i^{(a)} 1_{\zeta} E_j^{(b)} = \sum_{t \geq 0} \phi({a+b-\langle \check{\alpha}_i, \zeta\rangle \brack t}_i) E_i^{(a-t)} 1_{\zeta -(a+b-t)\alpha_i}F_i^{(b-t)}$;

\item $E_i 1_\zeta = 1_{\zeta+ \alpha_i}E_i$, $F_i 1_{\zeta} = 1_{\zeta - \alpha_i}F_i$;

\item $(u^+ 1_\zeta)(1_{\zeta'}u^-) = \delta_{\zeta, \zeta'}u^+ 1_\zeta u^-$ and $(u^- 1_\zeta)(1_{\zeta'}u^+) = \delta_{\zeta, \zeta'}u^- 1_\zeta u^+$;

\item $_R \Ud$ is a left module for $_R \mathbf U^+$ and $_R \mathbf U^-$.
\end{itemize}
\end{lemma}
\begin{proof}
This is shown, in slightly different notation, in \cite[31.1.3]{L93}
\end{proof}

\subsection{The Contraction Homomorphism}
It is now a simple consequence of the existence of the homomorphism $Fr'$ and the previous lemmas  that we have the following map, which we name the contracting homomorphism after the construction of Littelmann \cite{Li}. 

\begin{prop}
\label{contracthom}
Let $\phi\colon \A \to R$ be a homomorphism which factors through the natural map $\A \to \A_\ell$ for some positive integer $\ell$. Then (in the notation of Section \ref{modified}) there is an injective homomorphism $c\colon _R\Ud^* \to$ $_R\Ud$ given on generators by $e_i^{(n)}1_\lambda \mapsto E_i^{(nl_i)} 1_\lambda$ and $f_i^{(n)}1_\lambda \mapsto F_i^{(nl_i)} 1_\lambda$ where $\lambda \in X^* \subset X$.
\end{prop}
\begin{proof}
It suffices to prove the result for the ring $R = \mathcal A_\ell$. By lemma \ref{present} and the existence of the map $Fr'$ we must only check that the appropriate relations are satisfied. But then Lemma \ref{binom} establishes the result. The fact that the map $c$ is an injection follows from the fact that it is a right inverse to Lusztig's quantum Frobenius map $Fr\colon _R \Ud \to $ $_R \Ud^*$.
\end{proof}

\begin{remark}
The reader will presumably notice that, apart from the statement of the proposition, very little in this section is new. The author imagines that the proposition was known to both Lusztig and Littelmann. Indeed the existence of the splitting map $c$ along with the equivalence of categories between unital $\Ud$ modules and weight modules for $\U$ gives an alternative proof of Theorem $1$ of \cite{Li} (though of course the technical ingredients are essentially the same). The inclusion $X^* \subset X$ ''explains'' why one gets a contracting action of $_R \U^*$ only on a \textit{subspace} of a $_R \Ud$-module. 
\end{remark}

\section{Generalized $q$-Schur algebras}
\label{gqschursection}

\subsection{Saturated sets and $\Ud_P$}
We assume in this section that our root datum $(X,Y)$ is of finite type (\textit{i.e.} the Cartan matrix $(a_{ij})$ is positive definite). With this assumption, we may talk about the set of dominant weights $X^+ = \{ \lambda \in X: \langle \check{\alpha}_i, \lambda\rangle \geq 0\}$. We will show that the map $c$ is descends to a map on generalized $q$-Schur algebras. We first need to recall the definition of a generalized $q$-Schur algebra, following Lusztig \cite[Chapter 29]{L93}.

For $\lambda \in X^+$ we define the Weyl module of highest weight $\lambda$ to be:
\[
\Lambda_\lambda = \Ud1_\lambda/\big( \sum_i \Ud E_i 1_\lambda + \sum_i  \Ud F_i^{(\langle \check{\alpha}_i \lambda \rangle+1)} 1_\lambda \big).
\]
It is shown in \cite{L93} that the nonzero images of the elements $\{b^-1_\lambda: b \in \mathbf B\}$
form a canonical basis $\mathbf B(\lambda)$ of $\Lambda_\lambda$. (By abuse of notation, $\mathbf B(\lambda)$ is thought of both as a basis of $\Lambda_\lambda$ and as a subset of $\mathbf B$.) Taking the $\A$-span of $\mathbf B(\lambda)$ inside $\Lambda_\lambda$ we get a natural integral form $_\A \Lambda_\lambda$, and so a specialization for any $\A$-algebra $R$ which we denote $_R \Lambda_\lambda$. Over the field $\mathbb Q(v)$ the Weyl modules are irreducible, and the category of finite dimensional representations of $\Ud$ (more precisely, the category of finite dimensional representations of type 1) is semisimple.  

Recall that $X$ has a partial order $\leq$ such that $\lambda \leq \mu$ if $\mu - \lambda = \sum_i n_i \alpha_i$ where the $n_i$ are nonnegative integers. Using the Weyl modules, Lusztig introduced a decreasing (with respect to $\leq$) $X^+$-filtration $(\Ud[\geq \lambda])_{\lambda \in X^+}$ on $\Ud$ which can be characterized as follows: Let $u \in \Ud$, then $u \in \Ud[\geq \lambda]$ if whenever $u$ acts by a nonzero endomorphism on $\Lambda_{\lambda_1}$ we have $\lambda \leq \lambda_1$ (see \cite[29.1.2]{L93} for more details). Lusztig also shows that each ideal in the filtration is spanned by a subset $\Bd[\geq \lambda]$ of the canonical basis $\dot{\mathbf B}$. Similarly the ideals 
\[
\Ud[>\lambda] = \sum_{\mu > \lambda} \Ud[\geq \mu],
\]
are based ideals. We set $\Bd[\lambda] = \Bd[\geq \lambda] - \Bd[> \lambda]$.It is easy to see that if $\phi_\lambda \colon \Ud[\geq \lambda] \to \text{End}(\Lambda_\lambda)$ then $\Ud[>\lambda]$ is the kernel of $\phi_\lambda$, and so $\Ud[\geq \lambda]/\Ud[ > \lambda]$ is finite dimensional, \textit{i.e.} $\Bd[\lambda]$ is a finite set.

\begin{definition}
\label{gqschur}
A subset $P \subset X^+$ is said to be \textit{saturated} if $X^+ \backslash P$ is finite and if for all $\lambda \in X$ such that there is a $\mu \in P$ with $\mu \leq \lambda$ we have $\lambda \in P$. Given a saturated set $P \subset X^+$ let $\Ud[P]$ be the sum of the ideals $\Ud[\geq \lambda]$ where $\lambda \in P$, and set $\Ud_P = \Ud/\Ud[P]$. The $\mathbb Q(v)$-algebra $\Ud_P$ is called a generalized $q$-Schur algebra.
\end{definition}

The algebras $\Ud_P$ are named ''generalized $q$-Schur algebras" in \cite{Doty}\footnote{This terminology was also used earlier by Lusztig in personal communication with the author.}. 
In that paper it is shown that they are quantizations of the generalized Schur algebras of Donkin \cite{Do86}. It follows from the above that $\Ud_P$ is a finite dimensional algebra, and moreover that it is a based ring in the sense of \cite{L95}. (It is also what is sometimes referred to as \textit{cellular} but this will not be important for us). Since the set $\Omega_P$ of weights $\lambda$ such that $1_\lambda \notin \Ud[P]$ is finite, the algebra $\Ud_P$, unlike $\Ud$, has a multiplicative identity given by $\sum_{\lambda \in \Omega_P} 1_\lambda$. Similarly it is not hard to see that $\Ud_P$ is also a quotient of the quantum group $\U$ -- there is a natural map where the element $K_\mu$ is sent to $\sum_{\lambda \in \Omega_P} v^{\langle \mu,\lambda \rangle} 1_\lambda$, and a simple Vandermonde determinant argument shows that this map is surjective (see, for example,  \cite[Lemma 4.3]{L95} for a closely related statement). 

 Let $_\A \Ud_P$ be the corresponding quotient $_\A \Ud /( \Ud[P] \cap {_\A \Ud})$ of the integral form $_\A \Ud$. The algebra $_\A \Ud_P$  is an integral form of $\Ud_P$, \textit{i.e.} the map $\mathbb Q(v) \otimes $ ${_\A \Ud_P} \to \Ud_P$ is an isomorphism, since the ideal $\Ud[P]$ is based.

\subsection{$Fr$, $c$, and the algebras $\Ud_P$}
\label{compatibilities}
We now wish to understand how the maps $Fr$ and $c$ interact with the filtration $(\Ud[\geq \lambda])_{\lambda \in X^+}$. The key is to show that the ideals $\Ud[P]$ have a very simple generating set. For a saturated set $P \subseteq X^+$, let $I_P$ be the (two-sided) ideal generated by the elements $\{1_\nu: \nu \in P\}$. The next result is a version of the main theorem in \cite{Doty}, which in turn is motivated by the ``refined Peter-Weyl theorem'' of Lusztig. 

\begin{prop}
\label{qschur}
 The ideals $I_P$ and $\Ud[P]$ of $\Ud$ coincide.
\end{prop}
\begin{proof}
Let $I_\lambda = I_P$ where $P$ is the saturation of $\lambda$, that is, $P = \{ \mu \in X^+: \lambda \leq \mu\}$. It is clearly sufficient to show that $\Ud[\geq \lambda] = I_\lambda$, and moreover, it is easy to see that $I_\lambda \subseteq \Ud[\geq \lambda]$. It is shown in \cite[Proposition 4.4]{L95} that 
for any element $b$ of $\Bd[\lambda]$ there are $b_1,b_2 \in \B(\lambda)$ (where here we think of $\mathbf B(\lambda)$ as a subset of $\mathbf B$) such that
\[
b_1^-1_\lambda \sigma(b_2)^+ \equiv b \mod \Ud[> \lambda].
\]
(Here, $\sigma\colon \mathbf f \to \mathbf f$ is the antiautomorphism which is the identity on the generators of $\mathbf f$).
Thus it follows that $I_{\lambda} + \Ud[> \lambda] =\Ud[\geq \lambda]$. Now $\Ud[> \lambda] = \sum_{\mu > \lambda} \Ud[\geq \mu]$, and so since we have $I_\mu \subseteq I_\lambda$ and $\Ud[\geq \mu] \subseteq \Ud[\geq \lambda]$ for any $\mu$ with $\mu \geq \lambda$, we see that   in fact for any $\mu > \lambda$ we have 
\begin{equation}
\label{filtrationeq}
I_\lambda + \Ud[\geq \mu] = \Ud[\geq \lambda].
\end{equation}
Now suppose that $u \in \Ud[\geq \lambda]$. Since $\Ud[\geq \lambda]$ is a based ideal, may write it as a linear combination $u = c_1 b_1 + c_2 b_2 + \ldots c_k b_k$ for $b_i \in \Bd[\geq \lambda]$, $c_i \in \mathbb Q(v)$. The filtration $(\Ud[\geq \mu])_{\mu \in X^+}$ has intersection the zero ideal, so there exists a $\mu \in X^+$ such that $b_i \notin \Ud[\geq \mu]$ for all $i, 1 \leq i\leq k$. The set $\Bd[\geq \lambda] - \Bd[\geq \mu]$ spans a vector space complement $V[\lambda,\mu]$ to $\Ud[\geq \mu]$ in $\Ud[\geq \lambda]$, so that $\Ud[\geq \lambda] = V[\lambda,\mu] \oplus \Ud[\geq \mu]$. By Equation \ref{filtrationeq} we may write $u = v+w$ where $v \in I_\lambda$ and $w \in \Ud[\geq \mu]$, and moreover, clearly we can assume that $v \in V[\lambda,\mu]$. But since $u \in V[\lambda,\mu]$ by our choice of $\mu$, it follows immediately that $w=0$, that is, $u=v \in I_\lambda$ as required.
\end{proof}

\begin{remark}
Of course, Proposition 4.4 of \cite{L95} is much finer than we actually need. The results of \cite{L95} on the cell structure of the algebra $\Ud$ also motivated the statement of the proposition --- the theory of cells for $\Ud$ matches precisely with that of the Weyl modules, so the elements $1_\lambda$ for $\lambda \in X^+$ are in bijection with the two-sided cells, and given the elementary nature of the asymptotic algebra of each cell (a matrix algebra over $\mathbb Z$) the fact that the elements $\{1_\mu \in X^+ \colon \lambda \leq \mu \}$ generate $\Ud[\geq \lambda]$ is perhaps not surprising.
Note also that, although it is not necessary for the above proof, in fact $b^{\pm}1_\lambda \in \Bd$ for all $b \in \B$. 

Essentially the same result is contained, in slightly different notation, in \cite[Theorem 4.2]{Doty}, where it is established using the theory of finite dimensional algebras. The proof here has the advantage of generalizing to the affine case, as we will show in Section \ref{affine}.

\end{remark}

Now let $R$ be any $\A$-algebra, and define $_R \Ud_P$ to be the algebra $R \otimes {_\A \Ud_P}$. 
It follows from the previous lemma that these specializations have presentations analogous to Lemma \ref{present}. For $P$ a saturated set in $X^+$, set $P^* = X^* \cap P$. As in the previous section, we will be concerned with the case where the map $\A \to R$ factors through $\A_\ell$. In order to show that the map $c$ of section \ref{contract} descends to an embedding of generalized $q$-Schur algebras, we need only work over the ring $\A_\ell$.

\begin{theorem}
\label{embed}
Suppose that $R$ is an $\A_\ell$-algebra. The homomorphism $c$ induces an embedding $c_P\colon _R\Ud_{P^*}^* \to $ $_R\Ud_P$. Moreover, the quantum Frobenius map $Fr$ also induces a map $Fr_P \colon {_R\Ud_P} \to {_R\Ud^*_{P^*}}$.
\end{theorem}

\begin{proof}
As discussed above we may assume that $R=\A_\ell$. First note that from the definition of $I_P$ and $I_{P^*}^*$ it is clear that $c(I_{P^*}^*) \subseteq I_P$, hence $c$ induces a map $c_P \colon {_R\Ud^*_{P^*}} \to {_R\Ud_P}$. Similarly, if  $Fr\colon \Ud \to \Ud^*$ is the quantum Frobenius map, it follows immediately from the action of $Fr$ on generators that $Fr(I_P) = I_{P^*}^*$. It then also follows that $c_P$ is an embedding: suppose that $u \in \Ud^*$ is such that $c(u) \in I_P$, then we have $u = Fr(c(u)) \in I_{P^*}^*$ as required.
\end{proof}

\begin{remark}
Each of the algebras $\Ud_P$ is filtered by the images of the ideals $\Ud[\geq \lambda]$. The above proposition also immediately implies that the maps $c_P$ and $Fr_P$ are compatible with this filtration.
\end{remark}

\begin{remark}
Since the ideals $\Ud[P]$ are based, the results of this section suggest that there are compatibilities between the quantum Frobenius map and the canonical basis. The construction of the quantum Frobenius in \cite{McG06} gives a context in which it should be possible to investigate further compatibilities. I hope to return to this in a future paper.
\end{remark}

\section{$q$-Schur algebras}
\label{qschuralgebras}

\subsection{Geometric construction}
In this section, we relate our construction to \cite{BLM90} and show that we recover as a special case a construction of \cite{FM}. Thus we let $V$ be an $r$-dimensional vector space over a field $\mathbf k$, and $n$ a positive integer. Let $\mathcal F^n$ denote the space of $n$-step partial flags in $V$, that is 
\[
\mathcal F^n = \{ (0=F_0 \subseteq F_1 \subseteq \ldots \subseteq F_n = V): F_i \text{ a subspace of } V\}.
\]
Then the group $GL(V)$ acts transitively on the components $\mathcal F^n$, where the components are indexed by the set $\mathcal C_{n,r} = \{(a_1,a_2, \ldots, a_n) \in \mathbb N^n: \sum_{i=1}^n a_i = r\}$. The group $GL(V)$ also acts with finitely many orbits on $\mathcal F^n \times \mathcal F^n$. Here the orbits are indexed by the set $\Theta_r$ of $n\times n$ matrices $(a_{ij})$ with nonnegative integer entries such that $\sum_{i,j} a_{ij} = r$. If $(F,F') \in \mathcal F^n \times \mathcal F^n$ then the orbit it lies in is indexed by the matrix $(a_{ij})$ where 
\[
a_{ij} = \dim\biggl(\frac{F_i\cap F_j'}{(F_{i-1}\cap
F_j')+(F_i\cap F_{j-1}')}\biggr).
\]

Now suppose that $\mathbf k = \mathbb F_q$, a finite field with $q$ elements. Let $S_{\mathbf k}(n,r)$ denote the set of $\mathbb Z$-valued $GL(V)$-invariant functions on $\mathcal F^n \times \mathcal F^n$. Then $S_{\mathbf k}(n,r)$ is an algebra under convolution, and if we let $\mathbf 1 _A$ denote the indicator function for the $GL(V)$ orbit indexed by $A$ then $\{\mathbf 1_A: A \in \Theta_r\}$ is a $\mathbb Z$-basis of $S_{\mathbf k}(n,r)$. The structure constants of $S_\mathbf k(n,r)$ with respect to this basis are polynomial in $q$, hence we may define the $q$-Schur algebra $S_q(n,r)$ to be the $\mathbb Z[q]$-algebra such that $S_{\mathbf k}(n,r)$ is the specialization of  $S_q(n,r)$ at  $q = |\mathbf k|$ for any finite field $\mathbf k$. To connect this construction to the previous sections, we must extend scalars from $\mathbb Z[q]$ to $\mathbb Z[v,v^{-1}]$ by setting $q= v^2$. We will denote this extended algebra by $S_v(n,r)$. Let $(X,Y)$ be the root datum of type $SL_n$, so that $X = \mathbb Z^n/\mathbb Z(1,1,\ldots, 1)$ and $Y= \{(a_i) \in \mathbb Z^N: \sum_{i=1}^n a_i =0\}$, and for each $i \in \{1,2,\ldots, n-1\}$ we have $\alpha_i = (0,0,\ldots, 1,-1, \ldots 0) + \mathbb Z(1,1,\ldots,1)$ and $\check{\alpha}_i = (0,0,\ldots, 1,-1, \ldots 0)$. Let $\Ud$ the corresponding modified quantum group. Then the following is well known:

\begin{lemma}
The algebra $S_v(n,r)$ is isomorphic to $_\A \Ud_P$ where $P = \{\lambda \in X^+: \lambda \nleq  r\varpi_1\}$, and $\varpi_1$ is the first fundamental weight (the highest weight of the vector representation).
\end{lemma}
\begin{proof}
The proof that there is a homomorphism (given explicitly on the generators) $\psi_r \colon \Ud \to S_q(n,r)$ is implicit in \cite{BLM90}. A (simpler) proof that $\psi_r$ is a homomorphism in the affine case, which also applies in the finite type case  is given in \cite[7.7]{L99}.
The fact that the homomorphism is a surjection follows immediately from \cite[Proposition 3.9]{BLM90}. Showing that the kernel of $\psi_r$ is $I_P$, where $P$ is as in the statement of the Lemma, is essentially the $q$-analogue of Schur-Weyl duality (first established by Jimbo), an account of which is given in \cite{GL93}.
\end{proof}

\subsection{Contraction and Frobenius on $S_q(n,r)$}
Let $S_\mathbb Z(n,r)$ denotes the algebra $S_v(n,r)$ specialized at $q=1$ (it is the integral form of the classical Schur algebra) and let $S_\ell(n, \ell r)$ denote the specialization $\A_\ell \otimes_\A S_v(n,r)$. 
Thus the previous section shows that the contracting homomorphism $c$ yields a family of embeddings $c_r\colon S_\mathbb Z(n,r) \to S_\ell(n, \ell r)$, with right inverses $F_r$. We wish to characterize the maps $F_r$ and $c_r$ in the context of \cite{BLM90}.

The partial ordering on orbits given by closure induces a partial ordering on the elements of $\Theta_r$ which we will denote by $\preceq$. Let $\langle A \rangle$ denote the element of $S_v(n,r)$ which specializes to $\mathbf 1_A \in S_{\mathbf k}(n,r)$ for any finite field $\mathbf k$. Note that the element $1_\lambda$ is sent to $\langle D \rangle$ where $D = (d_i\delta_{ij})$ and
\[
\lambda = (d_1,d_2,\ldots, d_n) \text{ mod } (1,1,\ldots, 1)
\]
if such $(d_i)$ exist, and to zero otherwise.
Set  $[A] = v^{-d_A}\langle A \rangle$ where 
\[
d_A = \sum_{i \geq k, j < l} a_{ij}a_{kl}. 
\]
 First we have the following corollary of Proposition \ref{embed} and the results of \cite{BLM90}.
\begin{lemma}
\label{characterization}
Let $c_r \colon S_{\mathbb Z}(n,r) \to S_{\ell}(n, \ell r)$ be the embedding and let $A =(a_{ij}) \in \Theta_r$ let $\ell A$ denote the matrix $(\ell a_{ij})$. We have 
\[
c_r([A]) = [\ell A] +\sum_{B \in \Theta_{\ell r}} c_{B,A} [B],
\]
where $B \prec \ell A$ and $c_{B,A} \in \A_\ell$. Moreover, the map $c_r$ is characterized by this property. 
\end{lemma}
\begin{proof}
Proposition $3.9$ in \cite{BLM90} shows that $S_q(n,r)$ is generated by basis elements $[A]$ which are the nonzero images of $\{ E_i^{(r)} 1_{\lambda}, F_i^{(r)}1_\lambda: i \in I, r \geq 0, \lambda \in X\}$ under the quotient map $\psi_r \colon \Ud \to S_q(n,r)$ (that is, the proposition shows that this map is indeed a surjection). A restatement of the proposition is that for any matrix of $A \in \Theta_{r}$ there is a monomial of the form 
\[
E_{i_1}^{( s_1)}E_{i_2}^{(s_2)}\ldots E_{i_k}^{(s_k)}F_{j_1}^{(t_1)}F_{j_2}^{(t_2)}\ldots F_{k_p}^{( t_p)} 1_\lambda,
\] 
($s_1, s_2, \ldots, s_k, t_1, t_2, \ldots, t_p \in \mathbb N$, $i_1, i_2, \ldots, i_k, j_1, j_2, \ldots, j_p \in I$, $\lambda \in X^*$), whose image under $\psi$ is of the form
\[
[A] + \sum_{B \prec A} a_{B}[B].
\]
The proof of the proposition immediately implies that the image of the monomial
\[
E_{i_1}^{(\ell s_1)}E_{i_2}^{(\ell s_2)}\ldots E_{i_k}^{(\ell s_k)}F_{j_1}^{(\ell t_1)}F_{j_2}^{(\ell t_2)}\ldots F_{k_p}^{(\ell t_p)} 1_\lambda,
\] 
in $S_q(n,r)$ has leading term $[\ell A]$ with respect to the partial ordering on $\Theta_{\ell r}$. The lemma then follows by induction on the partial order $\prec$. The fact that the map $c_r$ is characterized by this property follows from the fact that the elements $[A]$ which are the images of the elements $\{E_i^{(r)} 1_{\lambda}, F_i^{(r)}1_\lambda: i \in I, r \in \mathbb N, \lambda \in X\}$ correspond to closed orbits, and so are minimal in the partial ordering $\preceq$.
\end{proof}

Interestingly, while we can only characterize $c_r$ in this context, it is possible to describe the map $F_r$ completely:
\begin{lemma}
The map $F_r$ is given by 
\[
F_r([A]) =\left\{\begin{array}{cc}[B],  & \qquad \text{if } \exists B \in \Theta_{r} \text{ such that } A = \ell B, \\0, & \text{otherwise.}\end{array}\right.
\]
\end{lemma}
\begin{proof}
We do not have a proof of this using the techniques of this paper, however \cite{McG06} gives a simple geometric proof. 
\end{proof}

\subsection{The case $n=2$}
\label{n=2case}
The paper \cite[Theorem 4.2]{FM} constructs an embedding $S_p(2,r) \to S_p(2,pr)$ where $S_p(n,r)$ denotes the algebra $S_\mathbb Z(n,r)$ base-changed to a field of characteristic $p$. Since any such field is an $\A_p$-algebra in a natural way, we obtain such an embedding from $c_r$ by base change. In \cite{FM} their map is computed explicitly on the basis corresponding to $\{[A]: A \in \Theta_r\}$ (by abuse of notation, we use $[A]$ to denote the images of the basis elements $[A] \in S_q(n,r)$ under the specialization map). Indeed they show that the map defined by
\[
\left[\begin{array}{cc}a & b \\c & d\end{array}\right] \mapsto 
\left\{\begin{array}{cc} \left[\begin{array}{cc}pa & pb \\pc & pd\end{array}\right], & \text{if } b \text{ or } c = 0,\\
\sum_{\epsilon=0}^{p-1} \left[\begin{array}{cc}pa+\epsilon & pb-\epsilon \\pc-\epsilon & pd+ \epsilon\end{array}\right],& \text{otherwise}.
\end{array}\right.
\]
is an algebra homomorphism from $S_p(2,r) \to S_p(2,pr)$. 
\begin{cor}
\label{identification}
When $n=2$, the embedding of Fayers and Martin is the specialization of the map $c_r$.
\end{cor}
\begin{proof}
This is immediate from the characterization of the map $c_r$ in Lemma \ref{characterization}, and the definition of the Fayers-Martin map given above, once one observes that in the case $n=2$, the partial order $\preceq$ is given by $[A] \preceq [B]$ exactly when $r(A) =  r(B)$, $c(A)= c(B)$, and $a_{11} \geq b_{11}$.
\end{proof}

\begin{remark}
It follows from the Corollary that, after specialization, the map $c_r$ is given explicitly by the formulas of Fayers and Martin. Thus in this case we have a refinement of Lemma \ref{characterization}, showing that the specializations of the constants $c_{B,A}$ satisfy certain congruences. I do not know how to generalize this to give an explicit description of $c_r$. (For $n=2$ it is natural to conjecture that the constants $c_{B,A}$ are either $0$ or a power of $v$, but I do have a proof of this.)
\end{remark}

\begin{remark}
After the first version of this paper was written the author discovered the preprint version of \cite{F} which establishes the existence of a $q$-analog of the map in \cite{FM} for $S_q(2,r)$. The coincidence of the map constructed there and the one in this paper, along with the final conjecture of that paper follow from Theorem \ref{embed} and Lemma \ref{characterization} above in a similar fashion to Corollary \ref{identification}. Indeed one must only observe that if $[A][B]$ is a ``codeterminant'' in the language of \cite{F}, then $A$ and $B$ are minimal in the partial ordering $\preceq$.
\end{remark}

\section{An affine analog}
\label{affine}
\subsection{BLN algebras}
In this section we wish to discuss a generalization of the map $c_P$ to the affine case. The first issue here is that one needs a notion of ``affine $q$-Schur algebra''. No formal definition of such an algebra (outside of type $A$) has been given, but there are a number of candidates. A complexity which already presents itself in type $A$ is that there are at least two competing definitions, each of which first arose in the work of Ginzburg and Vasserot \cite{GV}, and Lusztig \cite{L99}. By a natural generalization of the the construction of $S_q(n,r)$ given in the last section, they define an algebra $\mathfrak A_{D,n,n}$ which receives a homomorphism from the quantum  group $\U(\widehat{\mathfrak {sl}}_n)$. However this map is no longer surjective, and so it is reasonable to consider its image $\U_{D,n,n}$ as an alternative candidate for the affine $q$-Schur algebra. A similar issue arises in the work of Nakajima \cite{N01}, where he studies convolution algebras of equivariant $K$-homology on quiver varieties. These algebras again receive a homomorphism from the quantum loop algebra which is not surjective. Hence in this context we are also presented with two reasonable definitions: either the geometric convolution algebras themselves, or the image of the homomorphism from the quantum loop algebra.

We propose an alternative definition, which is essentially equivalent to choosing ``smaller'' of the two options: the images of the quantum loop algebra (we discuss the relation to the algebras of Lusztig-Ginzburg-Vasserot and Nakajima  more precisely in \ref{relationtoconvolution} and \ref{quiver} below). Rather than attempting to justify this choice, we rather wish to emphasize that most of the known results for these algebras suggest that the behaviors of the two families of algebras differ less than one might first expect. Nakajima's result \cite{N01} that despite the nonsurjectivity one can still classify the simple finite-dimensional representations of quantum loop algebras via the simple representations of the convolution algebras, and the result of \cite{McG} that the cell structures of the two algebras which arise in the work of Lusztig \cite{L99}, \cite{L99a} and Ginzburg-Vasserot \cite{GV} are essentially identical are both examples of this phenomenon.

Our definition is the obvious generalization of Lusztig's algebras $\Ud/\Ud[P]$ in the context of the work of Beck and Nakajima \cite{BN} on cells in affine quantum groups at level zero. To say this more precisely, we need to introduce the root datum for an affine quantum group. Thus we suppose the symmetric matrix $(a_{ij})$ attached to the Cartan datum is of affine type (see \cite[2.1.3]{L93}), and that $i \cdot i = 2$ for all $i \in I$. We further assume that the homomorphism $\mathbb Z[I] \to Y$ given by $i \mapsto \check{\alpha}_i$ is surjective with one-dimensional kernel, while the map $\mathbb Z[I] \to X$ given by $i \to \alpha_i$ has finite index and one-dimensional kernel. This is the degenerate affine root datum, corresponding to the affine quantum group at level zero. Fix $i_0 \in I$ such that $I_0 = \{i \in I : i \neq i_0\}$ is a basis of $\mathbb Z[Y]$. The root datum $(X,Y,I-\{i_o\})$ is then the associated finite-type root datum. We therefore have a set of dominant weights $X^+$ for the finite-type algebra, that is, 
\[
X^+ = \{ \lambda \in X: \langle i, \lambda\rangle \in \mathbb Z_{\geq 0} \quad \forall i \in I- \{i_0\}\}.
\] 
The Weyl group $W$ associated to the Cartan datum (the affine Weyl group) acts on $X$ via a finite quotient. We write $\Uda$ for the associated modified quantum group.

For each dominant weight $\lambda \in X^+$ there is a module $V(\lambda)$ known as an extremal weight module, see for example \cite{K94} for their construction (where they are denoted $V(\lambda)^{max}$). Following \cite{BN}, define natural filtration $(\dot{\mathbb U}[ \geq \lambda])_{\lambda \in X^+}$ on the modified affine quantum group at level zero, indexed by $\lambda \in X^+$: we say $u \in \mathbb U$ is $\mathbb U[\geq \lambda]$ if $u$ acts by zero on every $V(\lambda')$ such that $\lambda' \ngeq \lambda$ (here $\geq$ is the usual dominance ordering on $X$ with respect to the finite-type root datum). Define $\mathbb U[> \lambda]$ similarly. Then we have the following theorem:

\begin{theorem}\cite{BN}
\label{cells}
The ideals $\mathbb U[\geq \lambda]$ are based. Moreover the subsets 
\[
\mathbb B[\lambda] = \{b \in \dot{\mathbb B}: b \in \mathbb U[\geq \lambda], b \notin \mathbb U[\geq \lambda'], \forall \lambda' > \lambda\}
\]
are precisely the two sided cells of $\dot{\mathbb U}$.
\end{theorem}

Note that it follows immediately from the definitions that $1_\lambda \in \dot{\mathbb B}[\lambda]$. We can thus make the following definition in direct analogy with Definition \ref{gqschur}:

\begin{definition}
Given a saturated set in $X^+$ let $\dot{\mathbb U}[P] = \sum_{\lambda \in P} \dot{\mathbb U}[\geq \lambda]$ and let  $\dot{\mathbb U}_P$ to be the quotient $\dot{\mathbb U}/\dot{\mathbb U}_P$. We call this $\mathbb Q(v)$-algebra the BLN-algebra attached to the set $P \subset X^+$.
\end{definition}

\begin{remark}
The term ``BLN-algebra'' (for Beck, Lusztig and Nakajima) is intended to avoid the rather cumbersome term "generalized affine $q$-Schur algebra''.
\end{remark}

The results of Section \ref{contract} apply in the context of this section, so that there is an embedding $c\colon _R\dot{\mathbb U}^* \to {_R\dot{\mathbb U}}$, and a Frobenius $Fr\colon  {_R\dot{\mathbb U}} \to {_R\dot{\mathbb U}^*}$ when $R$ is an $\mathcal A$-algebra such that the homomorphism $\mathcal A \to R$ factors through $\mathcal A_\ell$. Since the ideals $I_P$ are based (by \cite{BN}), they have a natural integral form, and so it makes sense to specialize them to obtain algebras $_R\dot{ \mathbb U}_P$. In particular we have algebras $_R\dot{\mathbb U}^*_P$ and  $ _R \dot{\mathbb U}_P$. We wish to show that the results of section \ref{gqschursection} have an affine analogue. 

\begin{prop}
\label{BLNalg}
Let $P$ be a saturated set in $X^+$ and let $I_P$ be the ideal of $\dot{\mathbb U}$ generated by $\{1_\mu: \mu \in P\}$. Then $I_P = \dot{\mathbb U}[P]$.
\end{prop}
\begin{proof}
As in Proposition \ref{qschur}, it is enough to show that $\dot{\mathbb U}[\geq \lambda] = I_\lambda$ where $I_\lambda = I_P$ with $P=\{\mu: \mu \geq \lambda\}$. This can be proved as in that case also once we have the results of \cite{BN} at hand, so we only record here the adjustments that need to be made to the  previous argument. Beck and Nakajima show (Proposition $6.27$) that if $b \in \dot{\mathbb U}[\lambda]$, then there are elements $b_1,b_2$ of $\dot{\mathbb B}$ such that 
\[
b = b_11_\lambda b_2^\sharp \text{ mod } \dot{\mathbb U}[>\lambda]
\]
where $\sharp$ is an antiautomorphism which preserves the filtration $(\dot{\mathbb U}[ \geq \lambda])_{\lambda \in X^+}$. It follows immediately that the ideal $\dot{\mathbb U}[\geq \lambda]$ can be written as $I_\lambda + \dot{\mathbb U}[> \lambda]$.  Just as in the finite-type case, we deduce that 
for any $\mu > \lambda$ we have $I_\lambda + \dot{\mathbb U}[\geq \mu] = \dot{\mathbb U}[\geq \lambda]$. The rest of the proof proceeds exactly as in Proposition \ref{qschur}, once it is observed that the filtration $(\dot{\mathbb U}[\geq \mu])_{\mu \in X^+}$ has intersection the zero ideal (this follows, for example, from the ``Peter-Weyl-type decomposition'' of the crystal structure of $\dot{\mathbb B}$ given in \cite[4.3]{BN}).
\end{proof}

Note that Proposition \ref{gqschur} can be rewritten to give a presentation of a generalized $q$-Schur algebra. Similarly, Proposition \ref{BLNalg} gives a presentation of BLN-algebras.
We may now conclude as before that the quantum Frobenius and the contraction map descend to BLN-algebras. We use the notation of Section \ref{compatibilities}.

\begin{cor}
\label{affinecompatibility} Let $R$ be an $\A_\ell$-algebra. The maps $c$ and $Fr$ are compatible with the ideals $I_{P^*}$ and $I_P$ in other words, they descend to maps
\[
c_P \colon _R \dot{\mathbb U}^*_{P^*} \to {_R\dot{\mathbb U}}_P, \qquad Fr_P \colon {_R\dot{\mathbb U}}_P \to {_R \dot{\mathbb U}^*_{P^*}},
\]
for any saturated set $P \subset X^+$. Moreover, $c_P$ is an embedding and $Fr_P$ is surjective. 
\end{cor}
\begin{proof}
This proceeds exactly as in the finite-type case.
\end{proof}

We end this section by relating the BLN-algebras more precisely to the various convolution algebras considered by Lusztig and Nakajima mentioned at the start of this section. 

\subsection{Relation to convolution algebras}
\label{relationtoconvolution}
We show first that in type $A$ our algebras $\mathbb U_{P}$ are generalizations of the algebras $\mathbf U_{D,n,n}$ studied in \cite{L99}, Recall the results of that paper (see also \cite{GV} and \cite{SV}): the geometry of $n$-step periodic lattices in $\mathbf k[\varepsilon, \varepsilon^{-1}]^{\oplus D}$ is used to construct a convolution algebra $\mathfrak A_{D,n,n}$ with a canonical basis, which receives a homomorphism $\psi_D \colon \dot{\mathbb U} \to \mathfrak A_{D,n,n}$, from the modified form of the level zero affine quantum group $\dot{\mathbb U}$ of type $A$. Moreover, it is shown that the canonical basis of $\mathfrak A_{D,n,n}$ has a natural subset (consisting of \textit{aperiodic} elements) which spans the image of the quantum group $\mathbf U_{D,n,n} = \psi_D(\dot{\mathbb U})$. It was conjectured in that paper, and subsequently shown by Schiffmann and Vasserot \cite{SV}, that the map $\psi_D$ is compatible with the canonical basis of $\dot{\mathbb U}$, so that if $b \in \dot{\mathbb B}$ then $\psi_D(b)$ is either zero or an aperiodic basis element of $\mathfrak A_{D,n,n}$. 

\begin{lemma}
The algebras $\mathbf U_{D,n,n}$ are precisely the BLN-algebras $\dot{\mathbb U}_{P}$ where
\[
P = \{\mu : \mu \nleq D\varpi_1\}
\]
\end{lemma}
\begin{proof}
From the compatibility of the map $\psi_D$ with the canonical basis, it follows that $\text{ker}(\psi_D)$ is a union of two-sided cells. From the definition of the map $\psi_D$ (see \cite[7.7]{L99}) it is immediate that for $\lambda \in X^+$ we have $1_\lambda \in \text{ker}(\psi_D)$ if and only if $\lambda \nleq D\varpi_1$. Now using either \cite{McG} or Theorem \ref{cells} we know that the two-sided cells are in bijection with $X^+$, where a two-sided cell $\mathbf c$ is attached to $\lambda$ precisely when $1_\lambda \in \mathbf c$. The lemma follows immediately.

\end{proof}

\begin{remark}
 When $D < n$ the presentation of $\mathbf U_{D,n,n}$ given by Proposition \ref{BLNalg} has been obtained independently in \cite[Theorem 2.6.1]{DG} (by completely different means). In that paper they work with the algebra $\mathfrak A_{D,n,n}$, but when $D<n$ the algebras $\mathfrak A_{D,n,n}$ and $\mathbf U_{D,n,n}$ coincide (indeed it is easy to check that when $D<n$ all the canonical basis elements of $\mathfrak A_{D,n,n}$ must be aperiodic).
 
Corollary \ref{affinecompatibility} can also be interpreted in the context of the construction of \cite{L99}. Because the quotient of the quantum group is a proper subalgebra of the ``affine $q$-Schur algebra'' $\mathfrak A_{D,n,n}$ of that paper, it is not immediately clear that the maps $c$ and $Fr$ are defined on all of $\mathfrak A_{D,n,n}$ . Using the techniques of \cite{McG06}, it should be possible to show that at least the map $Fr$ extends. 
\end{remark} 

\subsection{Relation to quiver varieties}
\label{quiver}
Assume now that the root datum $(X,Y,I-\{i_0\})$ is of (finite) type ADE. We describe the connection between our BLN-algebras and the geometry of quiver varieties. We mostly follow the notation of \cite{N04} (see  \cite{N01} for more details). For each $\lambda, \mu \in X^+$ with $\mu \geq \lambda$ there are quiver varieties $\mathfrak M(\mu, \lambda)$ and $\mathfrak M_0(\mu, \lambda)$, where the first is a smooth symplectic variety, and the second is an affine variety. There is a proper birational morphism $\pi\colon \mathfrak M(\mu,\lambda) \to \mathfrak M_0(\mu, \lambda)$. Moreover there is a stratification of $\mathfrak M_0(\mu,\lambda) = \sqcup_{\mu \leq \nu \leq \lambda} \mathfrak M_0^{\text{reg}}(\nu,\lambda)$ such that the $\mathfrak M_0(\mu, \lambda)$ form an inductive system which stabilizes at $\mu =0$  (in particular  $\mathfrak M_0^{\text{reg}}(\mu,\lambda)$ is a smooth open dense subvariety of $\mathfrak M_0(\mu,\lambda)$).  

Let $\mathfrak M_0(\lambda)$ denote the limit of this system, and let $\mathfrak M(\lambda)$ denote the disjoint union $\bigsqcup_{\mu \leq \lambda} \mathfrak M(\mu,\lambda)$. The variety $Z(\lambda)$ is defined to be the fiber product $\mathfrak M(\lambda)\times_{\mathfrak M_0(\lambda)} \mathfrak M(\lambda)$. If $\lambda = \sum_{i \in I_0}w_i \varpi_i$ where $\varpi_i$ are the fundamental weights, and we let $G_\lambda = \prod_{i \in I_0} \text{GL}_{w_i}(\mathbb C)$, then there is a natural action of $\mathbb C^* \times G_\lambda$-action on $\mathfrak M(\lambda, \mu)$ and $\mathfrak M_0(\lambda, \mu)$ which is compatible with the map $\pi$. Finally, there is a distinguished point (denoted $0$) in $\mathfrak M_0(\mu,\lambda)$ onto which the $\mathbb C^*$-action contracts $\mathfrak M_0(\mu,\lambda)$ and the fiber $\mathfrak L(\mu,\lambda)$ of $\pi$ over $0$ is a Lagrangian subvariety of $\mathfrak M(\mu,\lambda)$. 

In \cite[Section 4]{N04} Nakajima notes that one can rephrase the construction of \cite{N01} to give a map 
\[
\Phi_\lambda \colon _{\mathcal A}\dot{\mathbb U} \to K^{\mathbb C^* \times G_\lambda}(Z(\lambda))/\text{torsion}
\]
If $\mathfrak L(\lambda)$ denotes the disjoint union of the varieties $\mathfrak L(\lambda, \mu)$, then the convolution product makes $K^{\mathbb C^* \times G_\lambda}(\mathfrak L(\lambda))$ into a module for $K^{\mathbb C^* \times G_\lambda}(Z(\lambda))$, and Nakajima has shown \cite[4.2]{N04} that this module, thought of as a $_{\mathcal A}\dot{\mathbb U}$-module via the map $\Phi_\lambda$, is isomorphic to $_{\mathcal A}V(\lambda)$ (here one must use the correct normalization of the map $\Phi_\lambda$ -- see \cite[4.13]{N04a}), the $\mathcal A$-form of the extremal weight module. 
We have the following consequence of the work of Nakajima:

\begin{lemma}
Let $P$ denote the saturation of the weight $\lambda$, \textit{i.e} $P = \{\mu: \mu \nleq \lambda\}$. Then the map $\Phi_\lambda$ factors though $_{\mathcal A}\dot{\mathbb U}_P$.
\end{lemma}
\begin{proof}
It is enough to show that $\text{ker}(\Phi_\lambda) \subset {_\mathcal A\dot{\mathbb U}}[P]$.
Recall the discussion of \cite[4.2]{N04}: if
\[
x \in \mathfrak M_0^{\text{reg}}(\mu,\lambda)\subset \mathfrak M_0(\lambda)
\]
then it can be shown that a subgroup of $\mathbb C^*\times G_\lambda$ isomorphic to $\mathbb C^* \times G_\mu$ stabilizes $x$. Since the forgetful functor gives a homomorphism $K^{\mathbb C^* \times G_\lambda}(Z(\lambda)) \to K^{\mathbb C^*\times G_{\mu}}(Z(\lambda))$ we see that $K^{\mathbb C^* \times G_\mu}(\mathfrak M(\lambda)_x)$ is naturally a $_{\mathcal A}\dot{\mathbb U}$-module. It follows from \cite[Section 3]{N01} and the results above that 
\begin{equation}
\label{fibers}
K^{\mathbb C^* \times G_\mu}(\mathfrak M(\lambda)_x) \cong  {_{\mathcal A}V(\mu)}.
\end{equation}

It is shown in \cite[10(iii)]{N98} that $\mathfrak M_0^{\text{reg}}(\mu,\lambda) \neq \emptyset$ and if and only if $\mu$ is a weight of the irreducible highest weight representation $L(\lambda)$ of the finite type Lie algebra associated to the root datum $(Y,X,I_0)$. Since it is known that the weights of such a representation are saturated with respect to the dominance order, it follows that 
\begin{equation}
\label{nonemptycriterion}
\mathfrak M_0^{\text{reg}}(\mu, \lambda) \neq \emptyset \iff \mu \leq \lambda.
\end{equation}
Clearly (\ref{fibers}) and (\ref{nonemptycriterion}) imply the result.
\end{proof}

\begin{remark}
We conjecture that in fact the BLN-algebra $\dot{\mathbb U}_P$ of the previous Lemma is isomorphic to $\Phi_\lambda(_\mathcal A\dot{\mathbb U})$. The principal obstruction to proving this (at least for the author) is that it is not known that the $K$-homology of the varieties $Z(\lambda)$ have good properties (such as those shown for $\mathcal M(\lambda,\mu)$ and $\mathcal L(\lambda,\mu)$ in \cite[Section 7]{N01}). In type $A$ however, such results are available, and the conjecture can be verified in this case. Such a result combined with Corollary \ref{affinecompatibility} would suggest that one can interpret the quantum Frobenius map in the context of quiver varieties. We hope to return to this question in a future paper.
\end{remark}

\end{document}